\documentstyle[12pt]{article}

\font\twlgot =eufm10 scaled \magstep1
\font\egtgot =eufm8
\font\sevgot =eufm7

\font\twlmsb =msbm10 scaled \magstep1
\font\egtmsb =msbm8
\font\sevmsb =msbm7

\newfam\gotfam
\def\pgot{\fam\gotfam\twlgot}
\textfont\gotfam\twlgot
\scriptfont\gotfam\egtgot
\scriptscriptfont\gotfam\sevgot
\def\got{\protect\pgot}

\newfam\msbfam
\textfont\msbfam\twlmsb
\scriptfont\msbfam\egtmsb
\scriptscriptfont\msbfam\sevmsb
\def\Bbb{\protect\pBbb}
\def\pBbb{\relax\ifmmode\expandafter\Bb\else\typeout{You cann't use
Bbb in text mode}\fi}
\def\Bb #1{{\fam\msbfam\relax#1}}

\def\thebibliography#1{\bigskip\section*{\centerline{
\bf References}\\}\list
  {[\arabic{enumi}]}{\settowidth\labelwidth{#1}\leftmargin\labelwidth
    \advance\leftmargin\labelsep
    \usecounter{enumi}}
    \def\newblock{\hskip .11em plus .33em minus .07em}
    \sloppy\clubpenalty4000\widowpenalty4000
    \sfcode`\.=1000\relax}

\def\op#1{\mathop{{\it\fam0} #1}\limits}

\newcommand{\id}{{\rm Id\,}}

\newcommand{\hm}{{\rm Hom\,}}
\newcommand{\dif}{{\rm Diff\,}}

\newcommand{\bll}{\bullet}
\newcommand{\beq}{\begin{equation}}
\newcommand{\eeq}{\end{equation}}
\newcommand{\ben}{\begin{eqnarray}}
\newcommand{\een}{\end{eqnarray}}
\newcommand{\be}{\begin{eqnarray*}}
\newcommand{\ee}{\end{eqnarray*}}
\newcommand{\bea}{\begin{eqalph}}
\newcommand{\eea}{\end{eqalph}}

\newcommand{\cO}{{\cal O}}
\newcommand{\cA}{{\cal A}}
\newcommand{\cJ}{{\cal J}}

\newcommand{\gd}{{\got d}}

\newcommand{\cI}{{\cal I}}

\newcommand{\cZ}{{\cal Z}}

\newcommand{\cK}{{\cal K}}

\newcommand{\al}{\alpha}
\newcommand{\bt}{\beta}
\newcommand{\dl}{\delta}

\newcommand{\f}{\phi}

\newcommand{\Om}{\Omega}
\newcommand{\m}{\mu}

\newcommand{\g}{\gamma}
\newcommand{\G}{\Gamma}
\newcommand{\e}{\epsilon}
\newcommand{\ve}{\varepsilon}

\newcommand{\bb}{{\bf 1}}

\newcommand{\w}{\wedge}

\newcommand{\wh}{\widehat}
\newcommand{\ol}{\overline}

\newcommand{\dr}{\partial}

\newcommand{\ar}{\op\longrightarrow}

\newcommand{\ot}{\otimes}

\newenvironment{eqalph}{\stepcounter{equation}
\setcounter{equationa}{\value{equation}}
\setcounter{equation}{0}

\begin{eqnarray}}{\end{eqnarray}\setcounter{equation}{\value{equationa}}}

\newcounter{example}
\newcounter{remark}
\newcounter{theorem}
\newcounter{proposition}
\newcounter{lemma}
\newcounter{corollary}
\newcounter{definition}

\def\theremark{\arabic{remark}}

\def\thedefinition{\arabic{theorem}}

\newenvironment{proof}{{\it Proof.}}{
\medskip }

\newenvironment{theo}{\refstepcounter{theorem} \medskip{\bf
Theorem \thedefinition.}\it}{\medskip }
\newenvironment{prop}{\refstepcounter{theorem} \medskip{\bf
Proposition \thedefinition.}\it}{\medskip }

\newenvironment{defi}{\refstepcounter{theorem} \medskip{\bf
Definition \thedefinition.} \it}{\medskip }

\newcommand{\mar}[1]{}

\begin{document}
\hbox{}

{\parindent=0pt

\begin{center}

{\large \bf A dilemma of nonequivalent definitions of differential
operators in noncommutative geometry}
\bigskip

{\bf G. Sardanashvily}

\medskip

Department of Theoretical Physics, Moscow State University, 117234
Moscow, RUSSIA

e-mail: gennadi.sardanashvily@unicam.it

\end{center}

\bigskip

\bigskip

{\bf Abstract.} In contrast with differential operators on modules
over commutative and graded commutative rings, there is no
satisfactory notion of a differential operator in noncommutative
geometry.





}

{\parindent=18pt

\bigskip
\bigskip

\centerline{\bf 1. Introduction}

\bigskip

Let $\cK$ be a commutative ring and $\cA$ a $\cK$-ring (a unital
algebra with $\bb\neq 0$). Let $P$ and $Q$ be $\cA$-modules. We
address the notion of a $\cK$-linear $Q$-valued differential
operator on $P$. If a ring $\cA$ is commutative, there is the
conventional definition (Definition \ref{ws131}, Section 2) of a
differential operator on $\cA$-modules (Grothendieck \cite{grot},
Krasil'shchik et al \cite{kras} and Giachetta et al \cite{book05})
(see Koszul \cite{kosz}, Akman \cite{akman1} and Akman et al
\cite{akman2} for the equivalent ones). This definition is
straightforwardly extended to modules over graded commutative
rings (Definition \ref{ww2}, Section 3), i.e., supergeometry
(Giachetta et al \cite{book05}). However, there are several
nonequivalent definitions of differential operators in
noncommutative geometry (Dubois-Violette et al \cite{dublmp},
Borowiec \cite{bor97}, Lunts et al \cite{lunts} and Giachetta et
al \cite{book05}). We show that none of them is satisfactory at
all (Section 5).

The main problem is that derivations of a noncommutative ring
$\cA$ fail to form an $\cA$-module. Let us recall that, given a
$\cK$-algebra $\cA$ and an  $\cA$-bimodule $Q$, a $Q$-valued
derivation of $\cA$ is defined as a $\cK$-module morphism
$u:\cA\to Q$ which obeys the Leibniz rule
\mar{ws100}\beq
u(ab)=u(a)b +a u(b), \qquad a,b\in\cA. \label{ws100}
\eeq
However, there is another convention such that the Leibniz rule
takes the form $u(ab)=bu(a) +a u(b)$ (Lang \cite{lang}). A graded
commutative ring is a particular noncommutative ring. However, the
definition of a derivation of a graded commutative algebra
(Bartocci et al \cite{bart}) also differs from the expression
(\ref{ws100}). Therefore, supergeometry is not a particular case
of noncommutative geometry.

The Chevalley--Eilenberg differential calculus over a $\cK$-ring
$\cA$ provides an important example of differential operators on
modules over $\cA$. If $\cA$ is the $\Bbb R$-ring of smooth real
functions on a smooth manifold $X$, this differential calculus is
the de Rham complex of exterior forms on $X$.

Throughout the paper, a two-sided $\cA$-module $P$ is  called the
$(\cA-\cA)$-bimodule. A $(\cA-\cA)$-bimodule is said to be a
$\cA$-bimodule if it is central over the center $\cZ_\cA$ of
$\cA$, i.e., $ap=pa$, $a\in \cZ_\cA$, $p\in P$. If $\cA$ is
commutative, a $\cA$-bimodule is simply called an $\cA$-module.

\bigskip
\bigskip

\centerline{\bf 2. Differential operators on modules over a
commutative ring}

\bigskip

Let $\cA$ be a commutative $\cK$-ring. Given $\cA$-modules $P$ and
$Q$, the $\cK$-module $\hm_\cK (P,Q)$ of $\cK$-module
homomorphisms $\Phi:P\to Q$ is endowed with the two different
$\cA$-module structures
\mar{5.29}\beq
(a\Phi)(p):= a\Phi(p),  \qquad  (\Phi\bll a)(p) := \Phi (a
p),\qquad a\in \cA, \quad p\in P. \label{5.29}
\eeq
We further refer to the second one as the $\cA^\bll$-module
structure. Let us put
\mar{spr172}\beq
\dl_a\Phi:= a\Phi -\Phi\bll a, \qquad a\in\cA. \label{spr172}
\eeq

\begin{defi} \label{ws131} \mar{ws131}
An element $\Delta\in\hm_\cK(P,Q)$ is called a $k$-order
$Q$-valued differential operator on $P$ if
\be
\dl_{a_0}\circ\cdots\circ\dl_{a_k}\Delta=0
\ee
for any tuple of $k+1$ elements $a_0,\ldots,a_k$ of $\cA$.
\end{defi}

The set $\dif_k(P,Q)$ of these operators inherits the $\cA$- and
$\cA^\bll$-module structures (\ref{5.29}). In particular, zero
order differential operators obey the condition
\be
\dl_a \Delta(p)=a\Delta(p)-\Delta(ap)=0, \qquad a\in\cA, \qquad
p\in P,
\ee
and, consequently, they are $\cA$-module morphisms $P\to Q$. First
order differential operators $\Delta$ satisfy the condition
\mar{ws106}\beq
\dl_b\circ\dl_a\,\Delta(p)= ba\Delta(p) -b\Delta(ap)
-a\Delta(bp)+\Delta(abp) =0, \quad a,b\in\cA. \label{ws106}
\eeq

Let $P=\cA$. Any zero order $Q$-valued differential operator
$\Delta$ on $\cA$ is defined by its value $\Delta(\bb)$. There is
an isomorphism $\dif_0(\cA,Q)=Q$ via the association
\be
Q\ni q\mapsto \Delta_q\in \dif_0(\cA,Q),
\ee
where the zero order differential operator $\Delta_q$ is given by
the equality $\Delta_q(\bb)=q$. A first order $Q$-valued
differential operator $\Delta$ on $\cA$ fulfils the condition
\be
\Delta(ab)=b\Delta(a)+ a\Delta(b) -ba \Delta(\bb), \qquad
a,b\in\cA.
\ee
It is called a $Q$-valued derivation of $\cA$ if $\Delta(\bb)=0$,
i.e., the Leibniz rule (\ref{ws100}) holds. If $\Delta$ is a
derivation of $\cA$, then $a\Delta$ is well for any $a\in \cA$.
Hence, derivations of $\cA$ form a left $\cA$-module $\gd(\cA,Q)$.
Any first order differential operator on $\cA$ falls into the sum
\be
\Delta(a)= a\Delta(\bb) +[\Delta(a)-a\Delta(\bb)]
\ee
of the zero order differential operator $a\Delta(\bb)$ and the
derivation $\Delta(a)-a\Delta(\bb)$. Thus, there is the
$\cA$-module decomposition
\mar{spr156'}\beq
\dif_1(\cA,Q) = Q \oplus\gd(\cA,Q). \label{spr156'}
\eeq

Let $P=Q=\cA$. The module $\gd\cA$ of derivations of $\cA$ is also
a Lie algebra over the ring $\cK$ with respect to the Lie bracket
\be
u\circ u'-u'\circ u, \qquad u,u'\in \gd\cA.
\ee
Accordingly, the decomposition (\ref{spr156'}) takes the form
\be
\dif_1(\cA) = \cA \oplus\gd\cA.
\ee

For instance,  let $\cA=C^\infty(X)$ be an $\Bbb R$-ring of smooth
real functions on smooth manifold $X$. Let $P$ be a projective
$C^\infty(X)$-modules of finite rank. In accordance with the
well-known Serre--Swan theorem, $P$ is isomorphic to the module
$\G(E)$ of global sections of some vector bundle $E\to X$. In this
case, Definition \ref{ws131} restarts familiar theory of linear
differential operators on vector bundles over a manifold $X$. Let
$Y\to X$ be an arbitrary bundle over $X$. The theory of
(nonlinear) differential operators on $Y$ is conventionally
formulated in terms of jet manifolds $J^kY$ of sections of $Y\to
X$ (Palais \cite{palais}, Krasil'shchik et al \cite{kras}, Bryant
et al \cite{bry} and Giachetta et al \cite{cmp}). Accordingly,
there is the following equivalent reformulation of Definition
\ref{ws131} in terms of the jet modules $\cJ^k(P)$ of a module $P$
(Krasil'shchik et al \cite{kras} and Giachetta et al
\cite{book05}). If $P=\G(E)$, these jet modules are isomorphic to
modules of sections of the jet bundles $J^kE\to X$.

Given an $\cA$-module $P$, let us consider the tensor product
$\cA\otimes_\cK P$ of $\cK$-modules $\cA$ and $P$. We put
\be
\dl^b(a\otimes p):= (ba)\otimes p - a\otimes (b p), \qquad p\in P,
\qquad a,b\in\cA.
\ee
Let $\m^{k+1}$ be the submodule of $\cA\ot_\cK P$ generated by
elements $\dl^{b_0}\circ \cdots \circ\dl^{b_k}(a\otimes p)$.

\begin{defi} \label{ww4} \mar{ww4}
The $k$-order jet module $\cJ^k(P)$ of a module $P$ is  the
quotient of the $\cK$-module $\cA\otimes_\cK P$ by $\m^{k+1}$. The
symbol $c\ot_kp$ stands for its elements.
\end{defi}

In particular, the first order jet module $\cJ^1(P)$ consists of
elements $c\ot_1 p$ modulo the relations
\be
\dl^a\circ \dl^b(\bb\ot_1 p)= ab\otimes_1 p -b\otimes_1 (ap)
-a\otimes_1 (bp) +\bb\ot_1(abp) =0.
\ee
The $\cK$-module $\cJ^k(P)$ is endowed with the $\cA$- and
$\cA^\bll$-module structures
\be
b(a\ot_k p):= ba\ot_k p, \qquad b\bll(a\otimes_k p):= a\otimes_k
(bp).
\ee
There exists the module morphism
\mar{5.44}\beq
J^k: P\ni p\mapsto \bb\otimes_k p\in \cJ^k(P) \label{5.44}
\eeq
of the $\cA$-module $P$ to the $\cA^\bll$-module $\cJ^k(P)$ such
that, seen as an $\cA$-module, $\cJ^k(P)$ is generated by elements
$J^kp$, $p\in P$. The following holds (Krasil'shchik et al
\cite{kras}).

\begin{theo} \label{t6} \mar{t6}
Any $k$-order $Q$-valued differential operator $\Delta$ of on an
$\cA$-module $P$ uniquely factorizes
\be
\Delta: P\ar^{J^k} \cJ^k(P)\ar Q
\ee
through the morphism $J^k$ (\ref{5.44}) and some $\cA$-module
homomorphism ${\got f}^\Delta: \cJ^k(P)\to Q$.
\end{theo}

Theorem \ref{t6} shows that $\cJ^k(P)$ is the representative
object of the functor $Q\to \dif_1k(P,Q)$.  Its proof is based on
the fact that the morphism $J^k$ (\ref{5.44}) is a $k$-order
$\cJ^k(P)$-valued differential operator on $P$. Let us denote $J:
P\ni p\mapsto \bb\ot p\in \cA\ot P$. Then, for any ${\got
f}\in\hm_\cA (\cA\ot P,Q)$, we obtain
\mar{ww3}\beq
\dl_{b_0}\circ\cdots\circ\dl_{b_0}({\got
f}\circ J)(p)={\got f}(\dl^{b_0}\circ\cdots \dl^{b_0}\circ(\bb\ot
p)). \label{ww3}
\eeq
The correspondence $\Delta\mapsto {\got f}^\Delta$ defines an
$\cA$-module isomorphism
\mar{5.50}\beq
\dif_k(P,Q)=\hm_{\cA}(\cJ^k(P),Q). \label{5.50}
\eeq
This isomorphism leads to the above mentioned equivalent
reformulation of Definition \ref{ws131}.

\begin{defi} \label{ww1} \mar{ww1}
A $k$-order $Q$-valued differential operator on a module $P$ is an
$\cA$-module morphism of the $k$-order jet module $\cJ^k(P)$ of
$P$ to $Q$.
\end{defi}

As was mentioned above, the Chevalley--Eilenberg differential
calculus over a commutative ring $\cA$ provides an important
example of differential operators over $\cA$ (Giachetta et al
\cite{book05}). It is a subcomplex
\mar{ws102}\beq
0\to \cK\ar^{\rm in}\cA\ar^d \cO^1[\gd\cA]\ar^d
\cO^2[\gd\cA]\ar^d\cdots \label{ws102}
\eeq
of the Chevalley--Eilenberg complex $C^*[\gd\cA,\cA]$ of the Lie
$\cK$-algebra $\gd\cA$ (Fuks \cite{fuks}) where $\cO^k[\gd\cA]$
are modules of $\cA$-multilinear (but not all $\cK$-multilinear)
skew-symmetric maps
\mar{+840'}\beq
\f^k:\op\times^k \gd\cA\to \cA \label{+840'}
\eeq
and $d$ is the Chevalley--Eilenberg coboundary operator
\mar{+840}\ben
&& d\f^k(u_0,\ldots,u_k)=\op\sum^k_{i=0}(-1)^iu_i
(\f^k(u_0,\ldots,\wh{u_i},\ldots,u_k)) +\label{+840}\\
&& \qquad \op\sum_{i<j} (-1)^{i+j} \f^k([u_i,u_j],u_0,\ldots, \wh
u_i, \ldots, \wh u_j,\ldots,u_k), \nonumber
\een
where the caret $\wh{}$ denotes omission. A direct verification
shows that if $\f^k$ is an $\cA$-multilinear map, so is $d\f^k$.
The graded module $\cO^*[\gd\cA]$ is provided with the structure
of a graded $\cA$-algebra with respect to the product
\mar{ws103}\ben
&& \f\w\f'(u_1,...,u_{r+s})=
\op\sum_{i_1<\cdots<i_r;j_1<\cdots<j_s} {\rm sgn}^{i_1\cdots
i_rj_1\cdots j_s}_{1\cdots r+s} \f(u_{i_1},\ldots,
u_{i_r}) \f'(u_{j_1},\ldots,u_{j_s}), \label{ws103} \\
&& \f\in \cO^r[\gd\cA], \qquad \f'\in \cO^s[\gd\cA], \qquad u_k\in
\gd\cA, \nonumber
\een
where sgn$^{...}_{...}$ is the sign of a permutation. This product
obeys the relations
\mar{ws98,9}\ben
&& d(\f\w\f')=d(\f)\w\f' +(-1)^{|\f|}\f\w d(\f'),
\quad \f,\f'\in \cO^*[\gd\cA], \label{ws98}\\
&& \f\w \f' =(-1)^{|\f||\f'|}\f'\w \f. \label{ws99}
\een
Using the first one, one can easily justify that the
Chevalley--Eilenberg differential $d$ obeys the relations
(\ref{ws106}) and, thus, it is a first order
$\cO^{k+1}[\gd\cA]$-valued differential operator on
$\cO^k[\gd\cA]$, $k\in Bbb N$. In particular, we have
\mar{spr708}\beq
(d a)(u)=u(a), \qquad a\in\cA, \qquad u\in\gd\cA. \label{spr708}
\eeq
It follows that $d(\bb)=0$ and $d$ is a $\cO^1[\gd\cA]$-valued
derivation of $\cO^0[\gd\cA]=\cA$.

For instance, the Chevalley--Eilenberg differential calculus over
the $\Bbb R$-ring $C^\infty(X)$ is the de Rham complex of exterior
forms on a smooth manifold $X$ where $d$ is the exterior
differential.

\bigskip
\bigskip

\centerline{\bf 3. Differential operators in supergeometry}

\bigskip

As was mentioned above, Definition \ref{ws131} is
straightforwardly extended to linear differential operators on
modules over a graded commutative ring.

Unless otherwise stated, by a graded structure throughout this
Section is meant a $\Bbb Z_2$-graded structure, and the symbol
$[.]$ stands for the $\Bbb Z_2$-graded parity. Recall that an
algebra $\cA$ is called graded if it is endowed with a grading
automorphism $\g$ such that $\g^2=\id$. A graded algebra falls
into the direct sum $\cA=\cA_0\oplus \cA_1$ of two $\Bbb
Z$-modules $\cA_0$ and $\cA_1$ of even and odd elements such that
$\g(a)=(-1)^{[a]}a$, $a\in\cA$, and
\be
[aa']=([a]+[a']){\rm mod}\,2, \qquad a,a'\in\cA.
\ee
If $\cA$ is a graded ring, then $[\bb]=0$. A graded algebra $\cA$
is called graded commutative if
\be
aa'=(-1)^{[a][a']}a'a.
\ee
Of course, a commutative ring is a graded commutative ring where
$\cA=\cA_0$. Given a graded algebra $\cA$, a left graded
$\cA$-module $Q$ is a left $\cA$-module such that
\be
[aq]=([a]+[q]){\rm mod}\,2.
\ee
A graded module $Q$ is split into the direct sum $Q=Q_0\oplus Q_1$
of two $\cA_0$-modules $Q_0$ and $Q_1$ of even and odd elements.
Similarly, right graded modules are defined. If $\cA$ is a graded
commutative ring, a graded $\cA$-module can be provided with a
graded $\cA$-bimodule structure by letting
\be
qa:= (-1)^{[a][q]}aq, \qquad a\in\cA, \qquad q\in Q.
\ee

Let $\cK$ be a commutative ring and $\cA$ a graded commutative
$\cK$-ring. Let $P$ and $Q$ be graded $\cA$-modules. The graded
$\cK$-module $\hm_\cK (P,Q)$ of graded $\cK$-module homomorphisms
$\Phi:P\to Q$ can be endowed with the two graded $\cA$-module
structures
\mar{ws11}\beq
(a\Phi)(p):= a\Phi(p),  \qquad  (\Phi\bll a)(p) := \Phi (a
p),\qquad a\in \cA, \quad p\in P, \label{ws11}
\eeq
called $\cA$- and $\cA^\bll$-module structures, respectively. Let
us put
\mar{ws12}\beq
\dl_a\Phi:= a\Phi -(-1)^{[a][\Phi]}\Phi\bll a, \qquad a\in\cA.
\label{ws12}
\eeq
The following generalizes Definition \ref{ws131} (Giachetta et al
\cite{book05}).

\begin{defi} \label{ww2} \mar{ww2}
An element $\Delta\in\hm_\cK(P,Q)$ is said to be a $k$-order
$Q$-valued graded differential operator on $P$ if
\be
\dl_{a_0}\circ\cdots\circ\dl_{a_k}\Delta=0
\ee
for any tuple of $k+1$ elements $a_0,\ldots,a_k$ of $\cA$.
\end{defi}

The set $\dif_k(P,Q)$ of these operators inherits the graded
module structures (\ref{ws11}). In particular, zero order graded
differential operators obey the condition
\be
\dl_a \Delta(p)=a\Delta(p)-(-1)^{[a][\Delta]}\Delta(ap)=0, \qquad
a\in\cA, \qquad p\in P,
\ee
i.e., they are graded $\cA$-module morphisms $P\to Q$. First order
graded differential operator $\Delta$ satisfy the condition
\be
&& \dl_a\circ\dl_b\,\Delta(p)= ab\Delta(p)-
(-1)^{([b]+[\Delta])[a]}b\Delta(ap)-
(-1)^{[b][\Delta]}a\Delta(bp)+\\
&& \qquad (-1)^{[b][\Delta]+([\Delta]+[b])[a]}\Delta(abp) =0,
\qquad a,b\in\cA, \quad p\in P.
\ee

Let $P=\cA$. Any zero order $Q$-valued graded differential
operator $\Delta$ on $\cA$ is defined by its value $\Delta(\bb)$.
Then there is a graded $\cA$-module isomorphism $\dif_0(\cA,Q)=Q$
via the association
\be
Q\ni q\mapsto \Delta_q\in \dif_0(\cA,Q),
\ee
where $\Delta_q$ is given by the equality $\Delta_q(\bb)=q$. A
first order $Q$-valued graded differential operator $\Delta$ on
$\cA$ fulfils the condition
\be
\Delta(ab)= \Delta(a)b+ (-1)^{[a][\Delta]}a\Delta(b)
-(-1)^{([b]+[a])[\Delta]} ab \Delta(\bb), \qquad  a,b\in\cA.
\ee
It is called a $Q$-valued graded derivation of $\cA$ if
$\Delta(\bb)=0$, i.e., the graded Leibniz rule
\mar{ws10}\beq
\Delta(ab) = \Delta(a)b + (-1)^{[a][\Delta]}a\Delta(b), \qquad
a,b\in \cA, \label{ws10}
\eeq
holds (cf. the Leibniz rule (\ref{ws100})). One obtains at once
that any first order graded differential operator on $\cA$ falls
into the sum
\be
\Delta(a)= \Delta(\bb)a +[\Delta(a)-\Delta(\bb)a]
\ee
of a zero order graded differential operator $\Delta(\bb)a$ and a
graded derivation $\Delta(a)-\Delta(\bb)a$. If $\Delta$ is a
graded derivation of $\cA$, then $a\Delta$ is so for any $a\in
\cA$. Hence, graded derivations of $\cA$ constitute a graded
$\cA$-module $\gd(\cA,Q)$.

Let $P=Q=\cA$. The module $\gd\cA$ of graded derivation of $\cA$
is also a Lie $\cK$-superalgebra with respect to the superbracket
\be
[u,u']=u\circ u' - (-1)^{[u][u']}u'\circ u, \qquad u,u'\in \cA.
\ee
We have the graded $\cA$-module decomposition
\be
\dif_1(\cA) = \cA \oplus\gd\cA.
\ee

Since $\gd\cA$ is a Lie $\cK$-superalgebra, let us consider the
Chevalley--Eilenberg complex $C^*[\gd\cA;\cA]$ of $\gd\cA$ where
the graded commutative ring $\cA$ is a regarded as a
$\gd\cA$-module (Fuks \cite{fuks}). It reads
\mar{ws85}\beq
0\to \cA\ar^{\dl^0}C^1[\gd\cA;\cA]\ar^{\dl^1}\cdots
C^k[\gd\cA;\cA]\ar^{\dl^k}\cdots \label{ws85}
\eeq
where $C^k[\gd\cA;\cA]$ are $\gd\cA$-modules of $\cK$-linear
graded morphisms of the graded exterior products $\op\w^k \gd\cA$
of the $\cK$-module $\gd\cA$ to $\cA$. Let us bring homogeneous
elements of $\op\w^k \gd\cA$ into the form
\be
\ve_1\w\cdots\w\ve_r\w\e_{r+1}\w\cdots\w \e_k, \qquad
\ve_i\in\gd\cA_0, \quad \e_j\in\gd\cA_1.
\ee
Then the coboundary operators of the complex (\ref{ws85}) are
given by the expression
\be
&& \dl^{r+s-1}c(\ve_1\w\cdots\w\ve_r\w\e_1\w\cdots\w\e_s)=
\op\sum_{i=1}^r (-1)^{i-1}\ve_i
c(\ve_1\w\cdots\wh\ve_i\cdots\w\ve_r\w\e_1\w\cdots\e_s)+
\nonumber \\
&& \qquad\op\sum_{j=1}^s (-1)^r\ve_i
c(\ve_1\w\cdots\w\ve_r\w\e_1\w\cdots\wh\e_j\cdots\w\e_s) + \\
&& \qquad\op\sum_{1\leq i<j\leq r} (-1)^{i+j}
c([\ve_i,\ve_j]\w\ve_1\w\cdots\wh\ve_i\cdots\wh\ve_j
\cdots\w\ve_r\w\e_1\w\cdots\w\e_s)+\\
&&\qquad \op\sum_{1\leq i<j\leq s} c([\e_i,\e_j]\w\ve_1\w\cdots\w
\ve_r\w\e_1\w\cdots
\wh\e_i\cdots\wh\e_j\cdots\w\e_s)+\\
&& \qquad \op\sum_{1\leq i<r,1\leq j\leq s} (-1)^{i+r+1}
c([\ve_i,\e_j]\w\ve_1\w\cdots\wh\ve_i\cdots\w\ve_r\w
\e_1\w\cdots\wh\e_j\cdots\w\e_s).
\ee
The subcomplex $\cO^*[\gd\cA]$ of the complex (\ref{ws85}) of
$\cA$-linear morphisms is the above mentioned graded
Chevalley--Eilenberg differential calculus over a graded
commutative $\cK$-ring $\cA$. Its coboundary operators $\dl^k$ are
first order $\cO^*{k+1}[\gd\cA]$-valued graded differential
operators on $\cO^*k[\gd\cA]$ (Giachetta et al \cite{book05}).

One may hope that Definition \ref{ww1} is also extended to modules
over graded commutative rings. One has introduced jets of these
modules (Giachetta et al \cite{book05}), but the corresponding
generalization of Theorem \ref{t6} has not been studied.

\bigskip
\bigskip

\centerline{\bf 4. Noncommutative differential calculus}

\bigskip

Turning to differential operators over a noncommutative ring, let
us start with a few examples of noncommutative differential
calculus.

In a general setting, a $\Bbb Z$-graded algebra $\Om^*$ over a
commutative ring $\cK$ is defined as a direct sum
\be
\Om^*= \op\oplus_k \Om^k, \qquad k=0,1,\ldots,
\ee
of $\cK$-modules $\Om^k$, provided with an associative
multiplication law $\al\cdot\bt$, $\al,\bt\in \Om^*$, such that
$\al\cdot\bt\in \Om^{|\al|+|\bt|}$, where $|\al|$ denotes the
degree of an element $\al\in \Om^{|\al|}$. It follows that $\Om^0$
is a (noncommutative) $\cK$-algebra $\cA$,  and $\Om^*$ is an
$(\cA-\cA)$-algebra. A graded algebra $\Om^*$ is called a
differential graded algebra or a differential calculus over $\cA=
\Om^0$ if it is a cochain complex of $\cK$-modules
\mar{spr260}\beq
0\to \cK\ar\cA\ar^\dl\Om^1\ar^\dl\cdots\Om^k\ar^\dl\cdots
\label{spr260}
\eeq
with respect to a coboundary operator $\dl$ which obeys the graded
Leibniz rule
\be
\dl(\al\cdot\bt)=\dl\al\cdot\bt +(-1)^{|\al|}\al\cdot \dl\bt.
\ee
In particular, $\dl:\cA\to \Om^1$ is a $\Om^1$-valued derivation
of a $\cK$-algebra $\cA$. One considers the minimal differential
graded subalgebra $\Om^*\cA$ of a differential calculus $\Om^*$
which contains $\cA$. It is called the minimal differential
calculus over $\cA$. Seen as an $(\cA-\cA)$-algebra, $\Om^*\cA$ is
generated by the elements $\dl a$, $a\in \cA$, and consists of
monomials $\al=a_0\dl a_1\cdots \dl a_k$, $a_i\in \cA$, whose
product obeys the juxtaposition rule
\be
(a_0\dl a_1)\cdot (b_0\dl b_1)=a_0\dl (a_1b_0)\cdot \dl b_1-
a_0a_1\dl b_0\cdot \dl b_1.
\ee

Let us generalize the Chevalley--Eilenberg differential calculus
over a commutative ring in Section 2 to a noncommutative
$\cK$-ring $\cA$. For this purpose, let us consider derivations
$u\in\gd\cA$ of $\cA$. They obey the Leibniz rule (\ref{ws100})
The set of these derivations $\gd\cA$ is both a $\cZ_\cA$-bimodule
and a Lie $\cK$-algebra with respect to the Lie bracket $[u,u']$,
$u,u'\in \gd\cA$. It is readily observed that derivations preserve
the center $\cZ_\cA$ of $\cA$.

Let us consider the Chevalley--Eilenberg complex $C^*[\gd\cA,\cA]$
of the Lie algebra $\gd\cA$ with coefficients in the ring $\cA$,
regarded as a $\gd\cA$-module. This complex contains a subcomplex
$\cO^*[\gd\cA]$ of $\cZ_\cA$-multilinear skew-symmetric maps with
respect to the Chevalley--Eilenberg coboundary operator $d$
(\ref{+840}). Its terms $\cO^k[\gd\cA]$ are $\cA$-bimodules. In
particular,
\mar{708'}\beq
\cO^1[\gd\cA]=\hm_{\cZ_\cA}(\gd\cA,\cA). \label{708'}
\eeq
The graded module $\cO^*[\gd\cA]$ is provided with the product
(\ref{ws103}) which obeys the relation (\ref{ws98}) and makes
$\cO^*[\gd\cA]$ into a differential graded algebra. One can think
of its elements as being noncommutative generalization of exterior
forms on a manifold. However, it should be noted that, if $\cA$ is
not commutative, there is nothing like the graded commutativity
(\ref{ws99}) in general. The minimal Chevalley--Eilenberg
differential calculus $\cO^*\cA$ over $\cA$ consists of the
monomials $a_0 da_1\w\cdots \w da_k$, $a_i\in\cA$, whose product
$\w$ (\ref{ws103}) obeys the juxtaposition rule
\be
(a_0d a_1)\w (b_0d b_1)=a_0d (a_1b_0)\w d b_1- a_0a_1d b_0\w d
b_1, \qquad a_i,b_i\in\cA.
\ee
For instance, it follows from the product (\ref{ws103}) that, if
$a,a'\in\cZ_\cA$, then
\mar{w123}\beq
da\w da'=-da'\w da, \qquad ada'=(da')a. \label{w123}
\eeq

\begin{prop} \label{w130} \mar{w130}
There is the duality relation
\mar{ws130}\beq
\gd\cA=\hm_{\cA-\cA}(\cO^1\cA,\cA).\label{ws130}
\eeq
\end{prop}

\begin{proof}
It follows from the definition (\ref{+840}) of the
Chevalley--Eilenberg coboundary operator that
\mar{spr708'}\beq
(d a)(u)=u(a), \qquad a\in\cA, \qquad u\in\gd\cA. \label{spr708'}
\eeq
This equality yields the morphism
\be
\gd\cA\ni u\mapsto \f_u\in \hm_{\cA-\cA}(\cO^1\cA,\cA), \qquad
\f_u(da):=u(a), \qquad a\in \cA.
\ee
This morphism is a monomorphism  because the module $\cO^1\cA$ is
generated by elements $da$, $a\in\cA$. At the same time, any
element $\f \in \hm_{\cA-\cA}(\cO^1\cA,\cA)$ induces the
derivation $u_\f(a):=\f(da)$ of $\cA$. Thus, there is a morphism
\be
\hm_{\cA-\cA}(\cO^1\cA,\cA)\to\gd \cA,
\ee
which is a monomorphism since $\cO^1\cA$ is generated by elements
$da$, $a\in\cA$.
\end{proof}

A different differential calculus over a noncommutative ring is
often used in noncommutative geometry (Connes \cite{conn} and
Landi \cite{land}). Let us consider the tensor product
$\cA\op\ot_\cK\cA$ of $\cK$-modules. It is brought into an
$\cA$-bimodule with respect to the multiplication
\be
b(a\ot a')c:=(ba)\ot (a'c), \qquad a,a',b,c\in\cA.
\ee
Let us consider its submodule $\Om^1(\cA)$ generated by the
elements $\bb\ot a-a\ot\bb$, $a\in\cA$. It is readily observed
that
\mar{w110}\beq
d:\cA\ni a\mapsto \bb\ot a -a\ot\bb\in\Om^1(\cA) \label{w110}
\eeq
is a $\Om^1(\cA)$-valued derivation of $\cA$. Thus, $\Om^1(\cA)$
is an $\cA$-bimodule generated by the elements $da$, $a\in\cA$,
such that the relation
\mar{w265}\beq
(da)b=d(ab)-adb, \qquad a,b\in \cA, \label{w265}
\eeq
holds. Let us consider the tensor algebra $\Om^*(\cA)$ of the
$\cA$-bimodule $\Om^1(\cA)$. It consists of the monomials
$a_0da_1\cdots da_k, \qquad a_i\in\cA$, whose product obeys the
juxtaposition rule
\be
(a_0d a_1)(b_0d b_1)=a_0d (a_1b_0)d b_1- a_0a_1d b_0 b_1, \qquad
a_i,b_i\in \cA,
\ee
because of the relation (\ref{w265}). The operator $d$
(\ref{w110}) is extended to $\Om^*(\cA)$ by the law
\be
d(a_0da_1\cdots da_k):=da_0da_1\cdots da_k,
\ee
that makes $\Om^*(\cA)$ into a differential graded algebra.

Of course, $\Om^*(\cA)$ is a minimal differential calculus. One
calls it the universal differential calculus over $\cA$ because of
the following property (Landi \cite{land}). Let $P$ be an
$\cA$-bimodule. Any $P$-valued derivation $\Delta$ of $\cA$
factorizes as $\Delta={\got f}^\Delta\circ d$ through some
$(\cA-\cA)$-module homomorphism ${\got f}^\Delta:\Om^1(\cA)\to P$.
Moreover, let $\cA'$ be another $\cK$-algebra and $(\Om'^*,\dl')$
its differential calculus over a $\cK$-ring $\cA'$. Any
homomorphism $\cA\to \cA'$ is uniquely extended to a morphism of
differential graded algebras $\rho^*: \Om^*(\cA) \to \Om'^*$ such
that $\rho^{k+1}\circ d=\dl' \circ\rho^k$. Indeed, this morphism
factorizes through the morphism of $\Om^*(\cA)$ to the minimal
differential calculus in $\Om'^*$ which sends $da\to \dl'\rho(a)$.
Elements of the universal differential calculus $\Om^*(\cA)$ are
called universal forms. However, they can not be regarded as the
noncommutative generalization of exterior forms because, in
contrast with the Chevalley--Eilenberg differential calculus, the
monomials $da$, $a\in\cZ_\cA$, of the universal differential
calculus do not satisfy the relations (\ref{w123}).

It seems natural to regard derivations of a noncommutative
$\cK$-ring $\cA$ and the Chevalley--Eilenberg coboundary operator
$d$ (\ref{+840}) as particular differential operators in
noncommutative geometry.

\bigskip
\bigskip

\centerline{\bf 5. Differential operators in noncommutative
geometry}

\bigskip

As was mentioned above, Definition \ref{ws131} provides a standard
notion of differential operators on modules over a commutative
ring. However, there exist its different generalizations to
modules over a noncommutative ring.

Let $P$ and $Q$ be $\cA$-bimodules over a noncommutative
$\cK$-ring $\cA$. The $\cK$-module $\hm_\cK(P,Q)$ of $\cK$-linear
homomorphisms $\Phi:P\to Q$ can be provided with the left $\cA$-
and $\cA^\bll$-module structures (\ref{5.29}) and the similar
right module structures
\mar{ws105}\beq
(\Phi a)(p):=\Phi(p)a, \qquad (a\bll\Phi)(p):=\Phi(pa), \quad
a\in\cA, \qquad p\in\ P. \label{ws105}
\eeq
For the sake of convenience, we will refer to the module
structures (\ref{5.29}) and (\ref{ws105}) as the left and right
$\cA-\cA^\bll$ structures, respectively. Let us put
\mar{ws133}\beq
\ol\dl_a\Phi:=\Phi a-a\bll\Phi, \qquad a\in\cA, \qquad \Phi\in
\hm_\cK(P,Q). \label{ws133}
\eeq
It is readily observed that
\be
\dl_a\circ\ol\dl_b=\ol\dl_b\circ\dl_a, \qquad a,b\in\cA.
\ee

The left $\cA$-module homomorphisms $\Delta: P\to Q$ obey the
conditions $\dl_a\Delta=0$, for all $a\in\cA$ and, consequently,
they can be regarded as left zero order $Q$-valued differential
operators on $P$. Similarly, right zero order differential
operators are defined.

Utilizing the condition (\ref{ws106}) as a definition of a first
order differential operator in noncommutative geometry, one
however meets difficulties. If $P=\cA$ and $\Delta(\bb)=0$, the
condition (\ref{ws106}) does not lead to the Leibniz rule
(\ref{ws100}), i.e., derivations of the $\cK$-ring $\cA$ are not
first order differential operators. In order to overcome these
difficulties, one can replace the condition (\ref{ws106}) with the
following one (Dubois-Violette et al \cite{dublmp}).

\begin{defi} \label{ws120} \mar{ws120}
An element $\Delta\in \hm_\cK(P,Q)$ is called a first order
differential operator on a bimodule $P$ over a noncommutative ring
$\cA$ if it obeys the condition
\mar{ws114}\ben
&& \dl_a\circ\ol\dl_b\Delta=\ol\dl_b\circ\dl_a\Delta=0,
\qquad  a,b\in\cA, \nonumber \\
&& a\Delta(p)b -a\Delta(pb) -\Delta(ap)b +\Delta(apb)=0, \qquad
p\in P. \label{ws114}
\een
\end{defi}

First order $Q$-valued differential operators on $P$ make up a
$\cZ_\cA$-module $\dif_1(P,Q)$.

If $P$ is a commutative bimodule over a commutative ring $\cA$,
then $\dl_a=\ol\dl_a$ and Definition \ref{ws120} comes to
Definition \ref{ws131} for first order differential operators.

In particular, let $P=\cA$. Any left or right zero order
$Q$-valued differential operator $\Delta$ is uniquely defined by
its value $\Delta(\bb)$. As a consequence, there are left and
right $\cA$-module isomorphisms
\be
&& Q\ni q\mapsto \Delta^{\rm R}_q\in\dif_0^{\rm R}(\cA,Q), \qquad
\Delta^{\rm R}_q(a)=qa, \qquad a\in\cA,\\
&& Q\ni q\mapsto \Delta^{\rm L}_q\in\dif_0^{\rm L}(\cA,Q), \qquad
\Delta^{\rm L}_q(a)=aq.
\ee
A first order $Q$-valued differential operator $\Delta$ on $\cA$
fulfils the condition
\mar{ws110}\beq
\Delta(ab)=\Delta(a)b+a\Delta(b) -a\Delta(\bb)b. \label{ws110}
\eeq
It is a derivation of $\cA$ if $\Delta(\bb)=0$. One obtains at
once that any first order differential operator on $\cA$ is split
into the sums
\be
&& \Delta(a)=a\Delta(\bb) +[\Delta(a)-a\Delta(\bb)], \\
&& \Delta(a)=\Delta(\bb)a +[\Delta(a)-\Delta(\bb)a]
\ee
of the derivations $\Delta(a)-a\Delta(\bb)$ or
$\Delta(a)-\Delta(\bb)a$ and the left or right zero order
differential operators $a\Delta(\bb)$ and $\Delta(\bb)a$,
respectively. If $u$ is a $Q$-valued derivation of $\cA$, then
$au$ (\ref{5.29}) and $ua$ (\ref{ws105}) are so for any
$a\in\cZ_\cA$. Hence, $Q$-valued derivations of $\cA$ constitute a
$\cZ_\cA$-module $\gd(\cA,Q)$. There are two $\cZ_\cA$-module
decompositions
\be
&& \dif_1(\cA,Q)= \dif_0^{\rm L}(\cA,Q) \oplus \gd(\cA,Q), \\
&& \dif_1(\cA,Q)= \dif_0^{\rm R}(\cA,Q) \oplus \gd(\cA,Q).
\ee
They differ from each other in the inner derivations $a\mapsto
aq-qa$.

Let $\hm_\cA^{\rm R}(P,Q)$ and $\hm_\cA^{\rm L}(P,Q)$ be the
modules of right and left $\cA$-module homomorphisms of $P$ to
$Q$, respectively. They are provided with the left and right
$\cA-\cA^\bll$-module structures (\ref{5.29}) and (\ref{ws105}),
respectively.

\begin{prop} \label{ws113} \mar{ws113}
An element $\Delta\in\hm_\cK(P,Q)$ is a first order $Q$-valued
differential operator on $P$ in accordance with Definition
\ref{ws120} iff it obeys the condition
\mar{n21}\beq
\Delta(apb)=(\op\dr^\to a)(p)b +a\Delta(p)b + a(\op\dr^\leftarrow
b)(p), \qquad p\in P, \quad  a,b\in\cA,\label{n21}
\eeq
where $\op\dr^\to$ and $\op\dr^\leftarrow$ are $\hm_\cA^{\rm
R}(P,Q)$- and $\hm_\cA^{\rm L}(P,Q)$-valued derivations of $\cA$,
respectively. Namely,
\be
(\op\dr^\to a)(pb)=(\op\dr^\to a)(p)b, \qquad (\op\dr^\leftarrow
b)(ap) =a(\op\dr^\leftarrow b)(p).
\ee
\end{prop}

\begin{proof}
It is easily verified that, if $\Delta$ obeys the equalities
(\ref{n21}), it also satisfies the equalities (\ref{ws114}).
Conversely, let $\Delta$ be a first order $Q$-valued differential
operator on $P$ in accordance with Definition \ref{ws120}. One can
bring the condition (\ref{ws114}) into the form
\be
\Delta(apb)=[\Delta(ap)-a\Delta(p)]b +a\Delta(p)b
+a[\Delta(pb)-\Delta(p)b],
\ee
and introduce the derivations
\be
(\op\dr^\to a)(p):= \Delta(ap)-a\Delta(p), \qquad
(\op\dr^\leftarrow b)(p):= \Delta(pb)-\Delta(p)b.
\ee
\end{proof}

For instance, let $P$ be a differential calculus over a $\cK$-ring
$\cA$ provided with an associative multiplication $\circ$ and a
coboundary operator $d$. Then $d$ exemplifies a $P$-valued first
order differential operator on $P$ by Definition \ref{ws120}. It
obeys the condition (\ref{n21}) which reads
\be
d(apb)=(da\circ p)b+a(dp)b + a((-1)^{|p|}p\circ db).
\ee
For instance, let $P=\cO^*\cA$ be the Chevalley--Eilenberg
differential calculus over $\cA$. In view of the relations
(\ref{708'}) and (\ref{ws130}), one can think of derivations
$u\in\gd\cA$ as being vector fields in noncommutative geometry. A
problem is that $\gd\cA$ is not an $\cA$-module. One can overcome
this difficulty as follows (Borowiec \cite{bor97}).

Given a noncommutative $\cK$-ring $\cA$ and an $\cA$-bimodule $Q$,
let $d$ be a $Q$-valued derivation of $\cA$. One can think of $Q$
as being a first degree term of a differential calculus over
$\cA$. Let $Q^*_{\rm R}$ be the right $\cA$-dual of $Q$. It is an
$\cA$-bimodule:
\be
(bu)(q):= bu(q), \qquad (ub)(q):=u(bq), \qquad  b\in\cA, \qquad
q\in Q.
\ee
One can associate to each element $u\in Q^*_{\rm R}$ the
$\cK$-module morphism
\mar{ws121}\beq
\wh u:\cA\in a\mapsto u(da)\in\cA. \label{ws121}
\eeq
This morphism obeys the relations
\mar{ws123}\beq
\wh{(bu)}(a) =bu(da), \qquad \wh u(ba)=\wh u(b)a+\wh{(ub)}(a).
\label{ws123}
\eeq
One calls $(Q^*_{\rm R},u\mapsto\wh u )$ the $\cA$-right  Cartan
pair,  and regards $\wh u$ (\ref{ws121}) as an $\cA$-valued first
order differential operator on $\cA$ (Borowiec \cite{bor97}). Let
us note that $\wh u$ (\ref{ws121}) need not be a derivation of
$\cA$ and fails to satisfy Definition \ref{ws120}, unless $u$
belongs to the two-sided $\cA$-dual $Q^*\subset Q^*_{\rm R}$ of
$Q$. Morphisms $\wh u$ (\ref{ws121}) are called into play in order
to describe (left) vector fields in noncommutative geometry
(Borowiec \cite{bor97} and Jara et al \cite{jara}).

In particular, if $Q=\cO^1\cA$, then $au$ for any $u\in\gd\cA$ and
$a\in\cA$ is a left noncommutative vector field in accordance with
the relation (\ref{spr708}).

Similarly, the $\cA$-left Cartan pair is defined. For instance,
$ua$ for any $u\in\gd\cA$ and $a\in\cA$ is a right noncommutative
vector field.

If $\cA$-valued derivations $u_1,\ldots u_r$ of a noncommutative
$\cK$-ring $\cA$ or the above mentioned noncommutative vector
fields $\wh u_1,\ldots \wh u_r$ on $\cA$ are regarded as first
order differential operators on $\cA$, it seems natural to think
of their compositions $u_1\circ\cdots u_r$ or $\wh u_1\circ\cdots
\wh u_r$ as being particular higher order differential operators
on $\cA$. Let us turn to the general notion of a differential
operator on $\cA$-bimodules.

By analogy with Definition \ref{ws131}, one may try to generalize
Definition \ref{ws120} by means of the maps $\dl_a$ (\ref{spr172})
and $\ol\dl_a$ (\ref{ws133}). A problem lies in the fact that, if
$P=Q=\cA$, the compositions $\dl_a\circ\dl_b$ and
$\ol\dl_a\circ\ol\dl_b$ do not imply the Leibniz rule and, as a
consequence, compositions of derivations of $\cA$ fail to be
differential operators.

This problem can be solved if $P$ and $Q$ are regarded as left
$\cA$-modules (Lunts et al \cite{lunts}). Let us consider the
$\cK$-module $\hm_\cK (P,Q)$ provided with the left $\cA-\cA^\bll$
module structure (\ref{5.29}). We denote by $\cZ_0$ its center,
i.e., $\dl_a\Phi=0$ for all $\Phi\in\cZ_0$ and $a\in\cA$. Let
$\cI_0=\ol \cZ_0$ be the $\cA-\cA^\bll$ submodule of $\hm_\cK
(P,Q)$ generated by $\cZ_0$. Let us consider:

(i) the quotient $\hm_\cK (P,Q)/\cI_0$,

(ii) its center $\cZ_1$,

(iii) the $\cA-\cA^\bll$ submodule $\ol \cZ_1$ of $\hm_\cK
(P,Q)/\cI_0$ generated by $\cZ_1$,

(iv) the $\cA-\cA^\bll$ submodule $\cI_1$ of $\hm_\cK (P,Q)$ given
by the relation $\cI_1/\cI_0=\ol \cZ_1$.

\noindent Then we define the $\cA-\cA^\bll$ submodules $\cI_r$,
$r=2,\ldots$, of $\hm_\cK (P,Q)$ by induction as
$\cI_r/\cI_{r-1}=\ol \cZ_r$, where $\ol \cZ_r$ is the
$\cA-\cA^\bll$ module generated by the center $\cZ_r$ of the
quotient $\hm_\cK (P,Q)/\cI_{r-1}$.

\begin{defi} \label{ws135} \mar{ws135}
Elements of the submodule $\cI_r$ of $\hm_\cK (P,Q)$ are said to
be left $r$-order $Q$-valued  differential operators on an
$\cA$-bimodule $P$ (Lunts et al \cite{lunts}).
\end{defi}

\begin{prop} \label{ws137} \mar{ws137}
An element $\Delta\in \hm_\cK (P,Q)$ is a differential operator of
order $r$ in accordance with Definition \ref{ws135} iff it is a
finite sum
\mar{ws138}\beq
\Delta(p)=b_i\Phi^i(p) +\Delta_{r-1}(p), \qquad b_i\in\cA,
\label{ws138}
\eeq
where $\Delta_{r-1}$ and $\dl_a\Phi^i$ for all $a\in\cA$ are
$(r-1)$-order differential operators if $r>0$, and they vanish if
$r=0$.
\end{prop}

\begin{proof}
If $r=0$, the statement is a straightforward corollary of
Definition \ref{ws135}. Let $r>0$. The representatives $\Phi_r$ of
elements of $\cZ_r$ obey the relation
\mar{ws139}\beq
\dl_c\Phi_r= \Delta'_{r-1}, \qquad  c\in \cA, \label{ws139}
\eeq
where $\Delta'_{r-1}$ is an $(r-1)$-order differential operator.
Then representatives $\ol\Phi_r$ of elements of $\ol\cZ_r$ take
the form
\be
\ol\Phi_r(p)=\op\sum_i c'_i\Phi^i(c_ip) + \Delta''_{r-1}(p),
\qquad c_i,c'_i\in \cA,
\ee
where $\Phi^i$ satisfy the relation (\ref{ws139}) and
$\Delta''_{r-1}$ is an $(r-1)$-order differential operator. Due to
the relation (\ref{ws139}), we obtain
\mar{ws140}\beq
\ol\Phi_r(p)=b_i\Phi^i(p) + \Delta'''_{r-1}(p), \quad b_i=c_ic'_i,
\quad \Delta'''_{r-1}=-\op\sum_ic'_i\dl_{c_i}\Phi^i +
\Delta''_{r-1}. \label{ws140}
\eeq
Hence, elements of $\cI_r$ modulo elements of $\cI_{r-1}$ take the
form (\ref{ws140}), i.e., they are given by the expression
(\ref{ws138}). The converse is obvious.
\end{proof}

If $\cA$ is a commutative ring, Definition \ref{ws135} comes to
Definition \ref{ws131}. Indeed, the expression (\ref{ws138}) shows
that $\Delta\in \hm_\cK (P,Q)$ is an $r$-order differential
operator iff $\dl_a\Delta$ for all $a\in\cA$ is a differential
operator of order $r-1$.

\begin{prop} \label{ws143} \mar{ws143}
If $P$ and $Q$ are $\cA$-bimodules, the set $\cI_r$ of $r$-order
$Q$-valued differential operators on $P$ is provided with the left
and right $\cA-\cA^\bll$ module structures.
\end{prop}

\begin{proof}
This statement is obviously true for zero order differential
operators. Using the expression (\ref{ws138}), one can prove it
for higher order differential operators by induction.
\end{proof}

Let $P=Q=\cA$. Any zero order differential operator on $\cA$ in
accordance with Definition \ref{ws135} takes the form $a\mapsto
cac'$ for some $c,c'\in\cA$.

\begin{prop} \label{ws146} \mar{ws146}
Let $\Delta_1$ and $\Delta_2$ be $n$- and $m$-order $\cA$-valued
differential operators on $\cA$, respectively. Then their
composition $\Delta_1\circ\Delta_2$ is an $(n+m)$-order
differential operator.
\end{prop}

\begin{proof}
The statement is proved by induction as follows. If $n=0$ or
$m=0$, the statement issues from the fact that the set of
differential operators possesses both left and right
$\cA-\cA^\bll$ structures. Let us assume that $\Delta\circ\Delta'$
is a differential operator for any $k$-order differential
operators $\Delta$ and $s$-order differential operators $\Delta'$
such that $k+s<n+m$. Let us show that $\Delta_1\circ\Delta_2$ is a
differential operator of order $n+m$. Due to the expression
(\ref{ws138}), it suffices to prove this fact when $\dl_a\Delta_1$
and $\dl_a\Delta_2$ for any $a\in\cA$ are differential operators
of order $n-1$ and $m-1$, respectively. We have the equality
\be
&& \dl_a(\Delta_1\circ\Delta_2)(b)=
a(\Delta_1\circ\Delta_2)(b)-(\Delta_1\circ\Delta_2)(ab)= \\
&& \qquad \Delta_1(a\Delta_2(b))+(\dl_a\Delta_1\circ\Delta_2)(b)-
(\Delta_1\circ\Delta_2)(ab) = \\
&&\qquad (\Delta_1\circ\dl_a\Delta_2)(b)+
(\dl_a\Delta_1\circ\Delta_2)(b),
\ee
whose right-hand side, by assumption, is a differential operator
of order $n+m-1$.
\end{proof}

Any derivation $u\in\gd\cA$ of a $\cK$-ring $\cA$ is a first order
differential operator in accordance with Definition \ref{ws135}.
Indeed, it is readily observed that
\be
(\dl_au)(b)= au(b)-u(ab)=-u(a)b, \qquad b\in\cA,
\ee
is a zero order differential operator for all $a\in\cA$. The
compositions $au$, $u\bll a$ (\ref{5.29}), $ua$, $a\bll u$
(\ref{ws105}) for any $u\in\gd\cA$, $a\in\cA$ and the compositions
of derivations $u_1\circ\cdots\circ u_r$ are also differential
operators on $\cA$ in accordance with Definition \ref{ws135}.

At the same time, noncommutative vector fields do not satisfy
Definition \ref{ws135} in general. First order differential
operators by Definition \ref{ws120} also need not obey Definition
\ref{ws135}, unless $P=Q=\cA$.

By analogy with Definition \ref{ws135} and Proposition
\ref{ws137}, one can define differential operators on right
$\cA$-modules as follows.

\begin{defi} \label{ws151} \mar{ws151}
Let $P$ and $Q$ be seen as right $\cA$-modules over a
noncommutative $\cK$-ring $\cA$. An element
$\Delta\in\hm_\cK(P,Q)$ is said to be a right zero order
$Q$-valued differential operator on $P$ if it is a finite sum
$\Delta=\Phi^i b_i$, $b_i\in\cA$, where $\ol\dl_a\Phi^i=0$ for all
$a\in\cA$. An element $\Delta\in\hm_\cK(P,Q)$ is called a right
differential operator of order $r>0$ on $P$ if it is a finite sum
\mar{ws150}\beq
\Delta(p)=\Phi^i(p)b_i +\Delta_{r-1}(p), \qquad b_i\in\cA,
\label{ws150}
\eeq
where $\Delta_{r-1}$ and $\ol\dl_a\Phi^i$ for all $a\in\cA$ are
right $(r-1)$-order differential operators.
\end{defi}

Definition \ref{ws135} and Definition \ref{ws151} of left and
right differential operators on $\cA$-bimodules are not
equivalent, but one can combine them as follows.

\begin{defi} \label{ws152} \mar{ws152}
Let $P$ and $Q$ be bimodules over a noncommutative $\cK$-ring
$\cA$. An element $\Delta\in\hm_\cK(P,Q)$ is a two-sided zero
order $Q$-valued differential operator on $P$ if it is either a
left or right zero order differential operator. An element
$\Delta\in\hm_\cK(P,Q)$ is said to be a two-sided differential
operator of order $r>0$ on $P$ if it is brought both into the form
\be
\Delta=b_i\Phi^i +\Delta_{r-1},\qquad b_i\in\cA,
\ee
and
\be
\Delta=\ol\Phi^i\ol b_i +\ol\Delta_{r-1}, \qquad \ol b_i\in\cA,
\ee
where $\Delta_{r-1}$, $\ol\Delta_{r-1}$ and $\dl_a\Phi^i$,
$\ol\dl_a\ol\Phi^i$ for all $a\in\cA$
  are two-sided $(r-1)$-order differential operators.
\end{defi}

One can think of this definition as a generalization of Definition
\ref{ws120} to higher order differential operators.

It is readily observed that two-sided differential operators
described by Definition \ref{ws152} need not be left or right
differential operators, and {\it vice versa}. At the same time,
$\cA$-valued derivations of a $\cK$-ring $\cA$ and their
compositions obey Definition \ref{ws152}.

In conclusion, note that Definition \ref{ww1} is also generalized
to a certain class of first order differential operators on
modules over a noncommutative ring.

The notion of a jet is extended to a module $P$ over a
noncommutative ring $\cA$. One can follow Definition \ref{ww4} in
order to define the left jets of a two-sided module $P$ over a
noncommutative ring $\cA$. However, the relation (\ref{ww3}) fails
to hold if $k>0$. Therefore, no module of left jets is the
representative object of left differential operators. The notion
of a right jet of $P$ meets the similar problem. Given an
$\cA$-module $P$, let us consider the tensor product
$\cA\otimes_\cK P\ot_\cK\cA$ of $\cK$-modules $\cA$ and $P$. We
put
\be
&& \dl^b(a\otimes p\ot c):= (ba)\otimes p\ot c - a\otimes (bp)\ot c,\\
&& \ol\dl^b(a\otimes p\ot c):= a\otimes p\ot (cb) - a\otimes
(pb)\ot c, \qquad p\in P, \qquad a,b,c\in\cA.
\ee
Let us denote by $\m^1$ the two-sided $\cA$-submodule of
$\cA\ot_\cK P\ot_\cK\cA$ generated by elements of the type
$\ol\dl^c\circ\dl^b(\bb\otimes p\ot\bb)$. One can define the first
order two-sided jet module $\cJ^1(P)=\cA\otimes_\cK
P\ot_\cK\cA/\m^1$ of $P$. Let us denote
\be
J: P\ni p\mapsto \bb\ot p\ot\bb\in \cA\ot P\ot\cA.
\ee
Then the equality
\be
\dl_{b_0}\circ \cdots \circ\dl_{b_k}({\got f}\circ J)(p) ={\got
f}(\dl^{b_0}\circ \cdots \circ\dl^{b_k} (\bb\ot p\ot\bb))
\ee
holds for any ${\got f}\in\hm_{\cA-\cA} (\cA\ot P,Q)$. It follows
that $\cJ^1(P)$ is the representative object of the functor $Q\to
\ol\dif_1(P,Q)$ where the $\cK$-module $\ol\dif_1(P,Q)$ consists
of two-sided first order differential operators $\Delta$ by
Definition \ref{ws152} which obey the condition
$\ol\dl^c\circ\dl^b\Delta=0$ for all $c,b\in\cA$. They are the
first order differential operators by (Dubois-Violette et al
\cite{dublmp}).

}

\end{document}